# Inverse source problem for the space-time fractional parabolic equation on a metric star graph with an integral overdetermination condition

R. R. Ashurov[1,3,a], Z. A. Sobirov[1,2,b], A. A. Turemuratova[2,4,c]

[1]V.I. Romanovskiy Institute of Mathematics, Uzbekistan Academy of Science, Tashkent, Uzbekistan

[2]National University of Uzbekistan, Tashkent, Uzbekistan,

[3]School of Engineering, Central Asian University, 264, Milliy Bog St., 111221, Tashkent, Uzbekistan

[4]Branch of Russian Economic University named after G. V. Plekhanov in Tashkent, Tashkent, Uzbekistan

E-mail: [a]ashurovr@gmail.com

[b]z.sobirov@nuu.uz, [c]ariuxanturemuratova@gmail.com

**Abstract** – In the present paper, we investigate the initial-boundary value problem for fractional order parabolic equation on a metric star graph in Sobolev spaces. First, we prove the existence and uniqueness results of strong solutions which are proved with the classical functional method based on a priori estimates. Moreover, the inverse source problem with the integral overdetermination condition for space-time fractional derivatives in Sobolev spaces is first considered in the present paper. By transforming the inverse problem to the operator-based equation, we showed that the corresponding resolvent operator is well-defined.

## 1. INTRODUCTION

The branched thin structures and metric graphs are widely used as a model in theoretical research on applied problems. They are used in the study of many complex systems from physics, biology, ecology, sociology, economics and finance [1, 2]. It is known that the initial and initial-boundary value problems on metric graphs are also used for the theoretical study of diffusion processes in branched structures and networks [3, 4]. Wave processes in simplest decorated graph and infinite tree graph are investigated in [5, 6]. The star graph has importance as an elementary part of any connected metric graph, i. e. any complex graph can be considered as a combination of star graphs. Investigating initial-boundary value problems give evidence to understand the nature of scattering at a certain branching point of the graph.

Fractional calculus, a branch of mathematics, finds applications in numerous scientific disciplines. It has proven to be a valuable tool for describing various processes in fields ranging from science to engineering. Over the years, researchers have dedicated significant attention to studying fractional differential equations (FDEs) on metric graphs, as evidenced by the works cited in references [7-9].

One notable contribution in this area is the work by D. Idczak et al. [10], where they utilized the Riemann-Liouville fractional derivative to define and characterize fractional Sobolev spaces. By establishing embedding results, they demonstrated the applicability of the obtained results to Sturm-Liouville fractional equations expressed in terms of the composition of right and left Riemann-Liouville fractional derivatives.

Moreover, inverse problems are one of the most important mathematical problems in science and mathematics. The problem of finding a solution to an inverse problem in differential equations was studied by many researchers. In [11], the authors investigate the uniqueness and stability of the solution to the inverse problem of determining the order of the Caputo time-fractional derivative for a subdiffusion equation. R. Ashurov et al. [12], focus on an inverse problem related to determining the orders of systems of fractional pseudo-differential equations and they analyze the uniqueness and stability of the solution and propose numerical methods for solving the problem. In [13], authors considered the nonlocal boundary problem $d_t^\rho u(t) + Au(t) = f(t)$, $(0 < \rho < 1, 0 < t \leq T)$, $u(\xi) = \alpha u(0) + \varphi$, ($\alpha$ is a constant and $0 < \xi \leq T$), in an arbitrary separable Hilbert space $H$ with the strongly positive selfadjoint operator $A$. They proved existence and uniqueness theorems for solutions of the problems under consideration. Furthermore, they investigated the inverse problems of determining the right-hand side of the equation and the function $\varphi$ in the boundary conditions.

Various approaches to inverse problems on metric graphs have been developed by several researchers, including P. Kurasov, S. Avdonin, G. Leugering, and others [14-18]. In [19], V. Yurko studied the inverse spectral problem for Sturm–Liouville operators on arbitrary compact graphs with standard gluing conditions at internal vertices. Notice, that all of the above works are devoted to the integer order differential equations.

Also, we refer to some works with the integral overdetermination condition. In [20, 21], Kamynin used the integral overdetermination condition for the solution of the inverse problem for a degenerate parabolic equation. Turdiev and Durdiev's [22] research focuses on the inverse coefficient problem for a time-fractional wave equation with the integral type overdetermination condition.

In the present paper, we focus on studying the direct and inverse problems on a metric star graph with integral type overdetermination condition. We consider the star metric graph $\Gamma$ with $n$ bonds consists of a finite set of vertices $V = \{v_k\}_0^n$ and a finite set of edges $E = \{e_k\}_1^n$, where $e_k$ connect the vertices $V_0$ and $V_k, k = \overline{1,n}$ [23] (see Fig. 1.) . By assigning the interval $(0, l_k)$, to the bond $e_k, k = \overline{1,n}$ of the graph we define the coordinates $x_k$ on each of the bonds. Where the vertex of the graph has a coordinate $0$ in each bond. Further, for simplicity, we will use $x$ instead of $x_k$.

We are interested in the following space-time fractional parabolic equations

$$\partial_{0,t}^\alpha u_k(x,t) + \partial_{x,l_k}^\beta \left(\gamma_k(x) D_{0,x}^\beta u_k(x,t)\right) = f(t) g_k(x,t) + h_k(x,t),$$

$$x \in \bar{e}_k, \ t \in (0,T], \ k = \overline{1,n}, \tag{1.1}$$

with the initial conditions

$$u_k(x,0) = 0, \ x \in \bar{e}_k, \ k = \overline{1,n}, \tag{1.2}$$

the vertex conditions

$$I_{0,x}^{1-\beta} u_k(0,t) = I_{0,x}^{1-\beta} u_j(0,t), \ t \in (0,T], \ k \neq j, \ k,j = \overline{1,n} \tag{1.3}$$

$$\sum_{k=1}^{n}\left(\gamma_k(0)D_{0,x}^{\beta}u_k(0,t)\right)=0, \ t\in(0,T], \tag{1.4}$$

the boundary conditions

$$I_{0,x}^{1-\beta}u_k(l_k,t)=0, \ t\in(0,T], \ k=\overline{1,n}, \tag{1.5}$$

where $\partial_{0,t}^{\alpha}$ denotes the left Caputo fractional derivative of order $\alpha\in(0,1)$ concerning the time variable $t$, $D_{0,x}^{\beta}$ and $\partial_{x,l_k}^{\beta}$, $k=\overline{1,n}$ stand for the left Riemann-Liouville and the right Caputo fractional derivative of order $\beta\in(1/2,1)$, respectively, with respect to the spatial variable $x$, $I_{0,x}^{1-\beta}$ is the left Riemann-Liouville fractional integral of order $1-\beta$. The real-valued functions $\gamma_k$ are absolutely continuous on $[0,l_k]$ and $0<p_1\leq|\gamma_k|\leq p_2$, $g_k(x,t)$, $h_k(x,t)$, $k=\overline{1,n}$, are given functions, $f(t)$ is an unknown function. So, we consider a problem of finding the pair of functions $\{u(x,t),f(t)\}$. For this purpose, we define an additional condition, more clearly, an integral overdetermination condition in the following way

$$\sum_{k=1}^{n}\int_{0}^{l_k}\eta_k(x)u_k(x,t)dx=\psi(t), \ t\in(0,T], \tag{1.6}$$

where $\eta_k(x), k=\overline{1,n}$, and $\psi(t)$ are known functions.

So, we separate the problem (1.1)-(1.6) into two problems, i.e., we search for a solution to the original problem in the form $\{u,f\}=\{z,f\}+\{y,0\}$.

**Problem $D$** (Direct problem). We call the problem finding a solution of the following equation

$$\partial_{0,t}^{\alpha}y_k(x,t)+\partial_{x,l_k}^{\beta}\left(\gamma_k(x)D_{0,x}^{\beta}y_k(x,t)\right)=h_k(x,t), \ k=\overline{1,n},$$

on each edge $e_k$, with initial conditions (1.2), the vertex conditions (1.3)-(1.4), and the boundary conditions (1.5), *the direct problem*.

**Problem $I$** (Inverse problem). We present *the inverse problem* as follows: find the pair of functions $\{z(x,t),f(t)\}$ in the following equation

$$\partial_{0,t}^{\alpha}z_k(x,t)+\partial_{x,l_k}^{\beta}\left(\gamma_k(x)D_{0,x}^{\beta}z_k(x,t)\right)=f(t)g_k(x,t), \ k=\overline{1,n}. \tag{1.7}$$

on each edge $e_k$, with initial, vertex, and boundary conditions, same as in (1.2)-(1.5), and an overdetermination condition, which is given by

$$\sum_{k=1}^{n}\int_{0}^{l_k}\eta_k(x)z_k(x,t)dx=E(t):=\psi(t)-\sum_{k=1}^{n}\int_{0}^{l_k}\eta_k(x)y_k(x,t)dx, \ t\in(0,T]. \tag{1.8}$$

Direct problems for similar equation which consist first (integer) order time-derivative and space fractional operator, same as in our case, were studied by V. Mehandiratta et al. [24], G. Leugering et al. [25], and Pasquini Soh Fotsing [26]. They also investigated the optimal control problems for parabolic fractional Sturm-Liouville equations on a star graph. They solved the problem using the method of Galerkin and obtained a weak solution. In our case, we use the functional method [27] to show the uniqueness and existence of a strong solution to the direct problem imposing the same conditions for the right-hand side of the equation.

The inverse problem for the fractional parabolic equation (1.1) in Sobolev spaces, even in the case of a interval on space variable, which is particular case of our result, does not appear in the literature.

The rest of this paper is organized as follows. In section 2, the paper provides definitions and lemmas necessary to understand the subsequent sections. The unique solvability of the direct problem is shown in section 3. It presents the derivation of the unique solution to the direct problem. In section 4, the paper addresses the solvability of the inverse problem by transforming it into an operator-based problem. It demonstrates that the corresponding resolvent operator is well-defined. Finally, in section 5, it is presented its conclusion, summarizing the findings and contributions of the study.

## 2. PRELIMINARIES

Now, we present the necessary concepts for the study of the problem.

For a function, $u: \Gamma \to \mathbb{R}$, defined on the graph, we put $u|_{e_k} = u_k$. For the functions defined on the graph, we also use vector-type notations $u = (u_1, \ldots, u_n)$, $u_x = \left(\dfrac{\partial u_1}{\partial x}, \ldots, \dfrac{\partial u_n}{\partial x}\right)$, $u_{xx} = \left(\dfrac{\partial^2 u_1}{\partial x^2}, \ldots, \dfrac{\partial^2 u_n}{\partial x^2}\right)$, $\int_\Gamma u \, d\Gamma = \sum_{k=1}^{n} \int_0^{l_k} u_k \, dx$. For $u: \Gamma \to \mathbb{R}$, $v: \Gamma \to \mathbb{R}$, we put $uv = (u_1 v_1, u_2 v_2, \ldots, u_n v_n)$.

Let $G_\tau = \{(x,t): x \in \Gamma, t \in (0, \tau]\}$, $0 \leq \tau \leq T$.

**Definition 2.1.** [28] The left and right fractional integrals of order $0 < \alpha < 1$ for a function $f \in L_1(0,T)$ are, respectively, defined by

$$D_{0,t}^{-\alpha} f(t) = I_{0,t}^{\alpha} f(t) := \frac{1}{\Gamma(\alpha)} \int_0^t \frac{f(\tau)}{(t-\tau)^{1-\alpha}} d\tau,$$

$$D_{t,T}^{-\alpha} f(t) = I_{t,T}^{\alpha} f(t) := \frac{1}{\Gamma(\alpha)} \int_t^T \frac{f(\tau)}{(\tau-t)^{1-\alpha}} d\tau.$$

**Proposition 2.1.** [24] If $\alpha > 0$ and $1 \leq p \leq \infty$, then $I_{0,t}^{\alpha}$ and $I_{t,T}^{\alpha}$ are continuous from $L_p(0,T)$ into itself and

$$\left\| I_{0,t}^{\alpha} f \right\|_{L_p(0,T)} \leq \frac{T^\alpha}{\Gamma(\alpha+1)} \| f \|_{L_p(0,T)}, \quad \left\| I_{t,T}^{\alpha} f \right\|_{L_p(0,T)} \leq \frac{T^\alpha}{\Gamma(\alpha+1)} \| f \|_{L_p(0,T)},$$

for all $f \in L_p(0,T)$.

**Definition 2.2.** [28] The left and right Caputo fractional derivatives of order $0 < \alpha < 1$ for a function $f$ on $[0,T]$ are respectively, defined by

$$\partial_{0,t}^{\alpha} f(t) = \frac{1}{\Gamma(1-\alpha)} \int_0^t \frac{f'(\tau)}{(t-\tau)^{\alpha}} d\tau, \qquad \partial_{t,T}^{\alpha} f(t) = \frac{-1}{\Gamma(1-\alpha)} \int_t^T \frac{f'(\tau)}{(\tau-t)^{\alpha}} d\tau,$$

provided that, the integrals, in the right-hand sides of these expressions, are exist.

**Definition 2.3.** [28] The left and right Riemann-Liouville fractional derivatives of order $0 < \alpha < 1$ for a function $f$ on $[0,T]$ are respectively, defined by

$$D_{0,t}^{\alpha} f(t) = \frac{1}{\Gamma(1-\alpha)} \frac{d}{dt} \int_0^t \frac{f(\tau)}{(t-\tau)^{\alpha}} d\tau,$$

$$D_{t,T}^{\alpha}f(t) = \frac{-1}{\Gamma(1-\alpha)}\frac{d}{dt}\int_{t}^{T}\frac{f(\tau)}{(\tau-t)^{\alpha}}d\tau,$$

provided that, the integrals, in the right-hand sides of these expressions, are exist.

From the definitions above, we have the following

$$I_{0,t}^{\alpha}D_{0,t}^{\alpha}f(t) = f(t) - f(0), \quad I_{0,t}^{\alpha}\partial_{0,t}^{\alpha}f(t) = f(t).$$

We define the following operators on the graph

$$\partial_{+}^{\beta}u = \{\partial_{0,x}^{\beta}u_1, \partial_{0,x}^{\beta}u_2, \ldots, \partial_{0,x}^{\beta}u_n\}, \quad \partial_{-}^{\beta}u = \{\partial_{x,l_1}^{\beta}u_1, \partial_{x,l_2}^{\beta}u_2, \ldots, \partial_{x,l_n}^{\beta}u_n\},$$

which can be considered as a fractional order of the gradient operator. The operators $D_{\pm}^{\beta}$ can be defined in a similar way.

We introduce some functional spaces that are useful in solving the considered problem. The fractional derivatives defined above are pointwise derivatives. For further study, we need to give a definition of the generalized fractional derivative, which is well-defined in some subspace of a fractional order Sobolev space. Such derivative is defined in [29].

Following to [29], by $H^{\alpha}(0,T)$, $0 < \alpha < 1$, we denote fractional Sobolev - Slobodeskii space governed by the norm (see [29], [30])

$$\|u\|_{H^{\alpha}(0,T)} := \left(\|u\|_{L_2(0,T)}^2 + \int_0^T\int_0^T\frac{|u(t)-u(s)|^2}{(t-s)^{1+2\alpha}}dsdt\right)^{\frac{1}{2}}.$$

We put

$$_0H^{\alpha}(0,T) = \{u \in H^{\alpha}(0,T) : u(0) = 0\} \text{ for } \frac{1}{2} < \alpha < 1,$$

$$H_{\alpha}(0,T) = \begin{cases} H^{\alpha}(0,T), & 0 \leq \alpha < \frac{1}{2}, \\ \{v \in H^{\frac{1}{2}}(0,T) : \int_0^T\frac{|v(t)|^2}{t}dt < \infty\}, & \alpha = \frac{1}{2}, \\ _0H^{\alpha}(0,T), & \frac{1}{2} < \alpha \leq 1. \end{cases}$$

The space $H_{\alpha}(0,T)$ is a Banach space with the norm [29]

$$\|v\|_{H_{\alpha}(0,T)} = \begin{cases} \|v\|_{H^{\alpha}(0,T)}, & 0 < \alpha < 1, \alpha \neq \frac{1}{2}, \\ \left(\|v\|_{H^{\frac{1}{2}}(0,T)}^2 + \int_0^T\frac{|v(t)|^2}{t}dt\right)^{\frac{1}{2}}, & \alpha = \frac{1}{2}. \end{cases}$$

According to [29] the space $_0C^1[0,T] = \{v \in C^1[0,T] : v(0) = 0\}$ is dense in $H_{\alpha}(0,T)$. In the same work, it is shown that $I_{0,t}^{\alpha} : L_2(0,T) \to H_{\alpha}(0,T)$ is injective and surjective, and so, the weak fractional derivative can be defined as $\partial_{0,t}^{\alpha} = (I_{0,t}^{\alpha})^{-1}$.

Furthermore, the following norms are equivalent in $H_{\alpha}(0,T)$ (see [29])

$$\|\partial_{0,t}^{\alpha}v\|_{L_2(0,T)} \sim \|v\|_{H_{\alpha}(0,T)}.$$

The weak right fractional derivative $\partial_{t,T}^{\alpha}$ with the domain $_{\alpha}H(0,T)$, which is the closure of the space $^{0}C^{1}[0,T] = \{v \in C^{1}[0,T] : v(T) = 0\}$ with respect to the norm $\|u\|_{H^{\alpha}(0,T)}$, can be defined analogically.

We notice, that in the case $\frac{1}{2} < \alpha < 1$ for any $v(t) \in H^{\alpha}(0,T)$ the weak Caputo derivatives can be defined by the equalities $\partial_{0,t}^{\alpha} v(t) = \partial_{0,t}^{\alpha}(v(t) - v(0))$ and $\partial_{t,T}^{\alpha} v(t) = \partial_{t,T}^{\alpha}(v(t) - v(T))$ (see [29]).

Above we defined the generalized (weak) fractional derivative of Caputo type. Now we need to define the generalized (weak) Riemann-Liouville derivative.

First, we introduce the spaces for only one section of the graph $(0,l)$ and then generalize. Let $0 < \beta < 1$. We define $AC_0^{\beta,2}(0,l) = AC_0^{\beta,2}([0,l],\mathbb{R})$ as the set of all functions that have representation

$$f(x) = \frac{c_0}{\Gamma(\beta)} x^{\beta-1} + I_{0,x}^{\beta} \phi(x) \text{ for a.e. } x \in [0,l],$$

with some $\phi \in L_2(0,l)$, $c_0 \in \mathbb{R}$, and by $AC_l^{\beta,2}(0,l) = AC_l^{\beta,2}([0,l],\mathbb{R})$ we mean the set of all functions that have the representation

$$g(x) = \frac{d_0}{\Gamma(\beta)} (l-x)^{\beta-1} + I_{0,x}^{\beta} \psi(x) \text{ for a.e. } x \in [0,l],$$

with some $\psi \in L_2(0,l)$, $d_0 \in \mathbb{R}$. It is clear, that by using these representations one can define right and left generalized fractional derivatives as those equal to $\phi(x)$ and $\psi(x)$, respectively, while the constants $c_0$ and $d_0$ represent the traces of the corresponding functions at $x = 0$ and $x = l$, respectively. Furthermore, for the case $0 < \beta < 1/2$ if $f \in L_2(0,l)$ and $f \in AC_0^{\beta,2}(0,l)$, then from these representations it directly follows, that the trace of the function $c_0$ is zero, and, similarly, if $g \in L_2(0,l)$ and $g \in AC_l^{\beta,2}(0,l)$, then $d_0 = 0$. So, for this case, one could not define a non-zero trace of the function.

For this reason, for the next consideration, we put $1/2 < \beta < 1$, which is necessary to define vertex conditions, in which traces are not necessarily equal to zero.

**Remark.** [25] We have the following characterization:

$$D_{0,x}^{\beta} f \in L_2(0,l) \Leftrightarrow f \in AC_0^{\beta,2}(0,l),$$

$$D_{x,l}^{\beta} f \in L_2(0,l) \Leftrightarrow f \in AC_l^{\beta,2}(0,l),$$

We set [10], [25]

$$H_+^{\beta}(0,l) = AC_0^{\beta,2} \cap L_2(0,l),$$

$$H_-^{\beta}(0,l) = AC_l^{\beta,2} \cap L_2(0,l).$$

Then it follows from the definitions of $AC_0^{\beta,2}(0,l)$ and $AC_l^{\beta,2}(0,l)$ that,

$$\rho \in H_+^{\beta} \Leftrightarrow \rho \in L_2(0,l) \text{ and } D_{0,x}^{\beta} \rho \in L_2(0,l),$$

$$\rho \in H_-^{\beta} \Leftrightarrow \rho \in L_2(0,l) \text{ and } D_{x,l}^{\beta} \rho \in L_2(0,l).$$

Furthermore, the space $H_+^{\beta}(0,l)$, equipped with the scalar product

$$(\varphi,\psi)_{H^\beta_+(0,l)} = (\varphi,\psi)_{L_2(0,l)} + \left(D^\beta_{0,x}\varphi, D^\beta_{0,x}\psi\right)_{L_2(0,l)},$$

and the space $H^\beta_-(0,l)$, equipped with the scalar product

$$(\varphi,\psi)_{H^\beta_-(0,l)} = (\varphi,\psi)_{L_2(0,l)} + \left(D^\beta_{x,l}\varphi, D^\beta_{x,l}\psi\right)_{L_2(0,l)},$$

form Hilbert spaces (see [10]).

**Lemma 2.1.** [24] Let $f \in F^\alpha := \{f \in C[0,T] : D^\alpha_{0,t}f \in L_2(0,T)\}$ and $g \in Q^\alpha := \{g \in C[0,T] : \partial^\alpha_{t,T}g \in L_2(0,T)\}$, then the following holds:

$$\int_0^T D^\alpha_{0,t}f(t)g(t)dt = \int_0^T f(t)\partial^\alpha_{t,T}g(t)dt + \left[g(t)I^{1-\alpha}_{0,t}f(t)\right]_{t=0}^{t=T}.$$

We put

$$\Upsilon_k(0,l_k) = \left\{\varphi \in H^\beta_+(0,l_k) : \partial^\beta_{x,l_k}\left(\gamma_k(x)D^\beta_{0,x}\varphi\right) \in L_2(0,l_k)\right\}, \quad k = \overline{1,n}.$$

Clearly, $\Upsilon_k(0,l_k)$, $k = \overline{1,n}$, equipped with the norm

$$\|\varphi\|^2_{\Upsilon_k(0,l_k)} = \|\varphi\|^2_{H^\beta_+(0,l_k)} + \left\|\partial^\beta_{x,l_k}\left(\gamma_k(x)D^\beta_{0,x}\varphi\right)\right\|_{L_2(0,l_k)},$$

respectively, are Hilbert spaces. This can be easily obtained from the completeness of the space $L_2(0,l)$.

Moreover, we introduce the space

$$V(G_T) = \left\{\begin{array}{l} u \in L_2(G_T) : u(\cdot,t) \in H_\alpha(0,T), \partial^\alpha_{0,t}u \in L_2(G_T) \\ u \in L_2(0,T; \oplus^n_{k=1}\Upsilon_k(0,l_k)), I^{1-\beta}_{0,x}u_k(0,t) = I^{1-\beta}_{0,x}u_j(0,t), k \neq j, \\ \sum_{k=1}^n \left(\gamma_k(0)D^\beta_{0,x}u_k(0,t)\right) = 0, \ I^{1-\beta}_{0,x}u_k(l_k,t) = 0, \quad k,j = \overline{1,n} \end{array}\right\}.$$

According to trace embedding theorems for fractional order Sobolev spaces (see [31], [32]), the traces in the definition of this space are well-defined. Consequently, the space $V(G_T)$, equipped with the norm

$$\|u\|^2_{V(G_T)} = \|u\|^2_{L_2(G_T)} + \left\|\partial^\alpha_{0,t}u\right\|^2_{L_2(G_T)} + \left\|D^\beta_{0,x}u\right\|^2_{L_2(G_T)} + \sum_{k=1}^n \left\|\partial^\beta_{x,l_k}\left(\gamma_k(x)D^\beta_{0,x}u_k\right)\right\|^2_{L_2(G_T)},$$

is also a Hilbert space.

**Lemma 2.2** [33] Let $v \in C[0,l_k]$, $1/2 < \beta < 1$ and $y, \omega \in j_k(0,l_k)$, $k = \overline{1,n}$. Then

$$\int_0^{l_k} \partial^\beta_{x,l_k}\left(v(x)D^\beta_{0,x}y(x)\right)\omega(x)dx =$$

$$= \int_0^{l_k} y(x)\partial^\beta_{x,l_k}\left(v(x)D^\beta_{0,x}\omega(x)\right)dx + \left[I^{1-\beta}_{0,x}y(x)\left(v(x)D^\beta_{0,x}\omega(x)\right)\right]_{x=0}^{x=l_k}$$

$$- \left[\left(v(x)D^\beta_{0,x}y(x)\right)I^{1-\beta}_{0,x}\omega(x)\right]_{x=0}^{x=l_k}.$$

We notice that, according to [31], [32], for $\beta > 1/2$ the traces of the functions exist.

## 3. UNIQUE SOLVABILITY OF THE DIRECT PROBLEM

**Theorem 3.1.** Let $h \in L_2(G_T)$. Then *Problem D* is uniquely solvable in $V(G_T)$.

**Proof.** Following to [27] we solve the direct problem for $y$ by bringing it into the operator equation:

$$Ay = h. \qquad (3.1)$$

We define our operator $A$ as

$$Ay\big|_{e_k} = \partial_{0,t}^\alpha y_k(x,t) + \partial_{x,l_k}^\beta \left(\gamma_k(x) D_{0,x}^\beta y_k(x,t)\right), \quad k = \overline{1,n}, \qquad (3.2)$$

with the domain

$$D(A) = V(G_T), \qquad (3.3)$$

and the range $R(A) \subset L_2(G_T)$.

**Proposition 3.1.** Operator $A: D(A) \to L_2(G_T)$ is continuous.

**Proof.** Continuity of the operator $A$ follows from the following inequality

$$\|Ay\|_{L_2(G_T)} = \left\|\partial_{0,t}^\alpha y + \partial_-^\beta\left(\gamma D_+^\beta y\right)\right\|_{L_2(G_T)} \leq \left\|\partial_{0,t}^\alpha y\right\|_{L_2(G_T)} + \left\|\partial_-^\beta\left(\gamma D_+^\beta y\right)\right\|_{L_2(G_T)} \leq \|y\|_{V(G_T)}. \qquad (3.4)$$

**Proposition 3.2.** The inverse operator $A^{-1}: R(A) \to V(G_T)$ is well-defined and continuous.

**Proof.** Taking into account the inequality [34]

$$\int_\Gamma \omega \partial_{0,t}^\alpha \omega \, d\Gamma \geq \frac{1}{2} \partial_{0,t}^\alpha \int_\Gamma \omega^2 \, d\Gamma, \qquad (3.5)$$

we square $Ay$ and integrate over $\Gamma$:

$$\int_\Gamma (Ay)^2 \, d\Gamma = \int_\Gamma \left[\partial_{0,t}^\alpha y + \partial_-^\beta\left(\gamma D_+^\beta y\right)\right]^2 d\Gamma = \int_\Gamma \left[\left(\partial_{0,t}^\alpha y\right)^2 + \left(\partial_-^\beta\left(\gamma D_+^\beta y\right)\right)^2\right] d\Gamma +$$

$$+ 2 \cdot \int_\Gamma \partial_{0,t}^\alpha y \cdot \left(\partial_-^\beta\left(\gamma D_+^\beta y\right)\right) d\Gamma = \int_\Gamma \left[\left(\partial_{0,t}^\alpha y\right)^2 + \left(\partial_-^\beta\left(\gamma D_+^\beta y\right)\right)^2\right] d\Gamma +$$

$$+ 2 \cdot \int_\Gamma \gamma D_+^\beta y \cdot D_+^\beta \left(\partial_{0,t}^\alpha y\right) d\Gamma \geq \int_\Gamma \left[\left(\partial_{0,t}^\alpha y\right)^2 + \left(\partial_-^\beta\left(\gamma D_+^\beta y\right)\right)^2\right] d\Gamma + p_1 \cdot \partial_{0,t}^\alpha \int_\Gamma \left(D_+^\beta y\right)^2 d\Gamma.$$

Integrating the last inequality with respect to $t$ we get

$$\int_0^t d\tau \int_\Gamma \left[\left(\partial_{0,t}^\alpha y\right)^2 + \left(\partial_-^\beta\left(\gamma D_+^\beta y\right)\right)^2\right] d\Gamma + p_1 \cdot \partial_{0,t}^{\alpha-1} \int_\Gamma \left(D_+^\beta y\right)^2 d\Gamma \leq \int_0^t d\tau \int_\Gamma (Ay)^2 \, d\Gamma. \qquad (3.6)$$

Taking into account the homogeneous initial condition (1.2), from the estimate (3.6) we conclude, that $Ay = 0$ iff $y = 0$. So, the inverse operator exists.

According to the definition $\partial_{0,\tau}^\alpha = \left(I_{0,\tau}^\alpha\right)^{-1}$ (see [29]). So, using Proposition 2.1, we obtain the following analog of the Friedrichs's inequality for the case of the Caputo fractional derivative

$$\|u\|_{L_2(0,t)} = \left\|I_{0,\tau}^\alpha\left(\partial_{0,\tau}^\alpha u\right)\right\|_{L_2(0,t)} \leq \frac{t^\alpha}{\Gamma(\alpha+1)} \left\|\partial_{0,\tau}^\alpha u\right\|_{L_2(0,t)}.$$

The last inequality together with the estimate (3.6) shows, that the inverse operator $A^{-1}$ is continuous.

From the Proposition 3.1 and Proposition 3.2 we can conclude the following result.

**Corollary.** The range $R(A)$ of the operator $A$ is a closed linear subspace of $L_2(G_T)$.

Now to show the solvability of the direct problem it is sufficient to prove that there is no orthogonal complement to the range $R(A)$ of the operator $A$ in $L_2(G_T)$.

**Lemma 3.1.** If for some $\omega \in L_2(G_T)$ it holds $(Av, \omega) = 0$ for all $v \in D(A)$, then $\omega = 0$.

**Proof.** We must prove that if for some $\omega \in L_2(G_T)$

$$\int_0^t d\tau \int_\Gamma \left( \partial_{0,t}^\alpha v + \partial_-^\beta \left( \gamma D_+^\beta v \right) \right) \omega \, d\Gamma = 0, \tag{3.7}$$

for all $v \in D(A)$ then $\omega = 0$.

Let $Lv = \partial_-^\beta \left( \gamma D_+^\beta v \right)$, with the domain

$$D(L) = \{ w \in \bigoplus_{k=1}^n \Upsilon_k, \sum_{k=1}^n \left( \gamma_k(0) D_{0,x}^\beta w_k(0) \right) = 0, I_{0,x}^{1-\beta} w_k \big|_{x=l_k} = 0,$$

$$I_{0,x}^{1-\beta} w_k \big|_{x=0} = I_{0,x}^{1-\beta} w_j \big|_{x=0}, k \neq j, k, j = \overline{1,n} \}.$$

On the edge $e_k$, $k = 1, 2, \ldots, n$, of the graph this operator acts as

$$(Lv)\big|_{e_k} = \partial_{x,l_k}^\beta \left( \gamma_k(x) D_{0,x}^\beta v_k(x,t) \right).$$

We choose $v(x,t) = \partial_{t_1,t}^{-\alpha} L^{-1} \omega$ for $0 < t_1 < t < T$ and $v = 0$ for $0 < t < t_1$, where $L^{-1}$ – inverse operator of $L$.

Below we give how to define $w = L^{-1}\omega$.

Let

$$\partial_{x,l_k}^\beta \left( \gamma_k(x) D_{0,x}^\beta w_k(x) \right) = \omega(x), \quad k = \overline{1,n},$$

and $w = (w_1, \ldots, w_n)$ satisfies the following boundary and vertex conditions

$$I_{0,x}^{1-\beta} w_k(l_k) = 0, \sum_{k=1}^n \left( \gamma_k(0) D_{0,x}^\beta w_k(0) \right) = 0, I_{0,x}^{1-\beta} w_k(0) = I_{0,x}^{1-\beta} w_k(0), k \neq j, k, j = \overline{1,n}. \tag{3.8}$$

The solution to this equation is in the following form

$$w_k(x) = I_{0,x}^\beta \left( \frac{1}{\gamma(x)} I_{x,l_k}^\beta \omega(x) \right) + a_k I_{0,x}^\beta \left( \frac{1}{\gamma(x)} \right) + \frac{x^{\beta-1}}{\Gamma(\beta)} b.$$

Unknown coefficients $\mu = (b, a_1, a_2, \ldots, a_n)^t$ can be found using the boundary and vertex conditions (3.8). For that, the determinant of the matrix $P$ must differ from zero in the equation

$$P\mu = M, \tag{3.9}$$

where

$$P = \begin{pmatrix} 0 & 1 & 1 & \ldots & 1 \\ 1 & \int_0^{l_1} \frac{dx}{\gamma_1(x)} & 0 & \ldots & 0 \\ 1 & 0 & \int_0^{l_2} \frac{dx}{\gamma_2(x)} & \ldots & 0 \\ \vdots & \vdots & \vdots & \ldots & \vdots \\ 1 & 0 & 0 & \ldots & \int_0^{l_n} \frac{dx}{\gamma_n(x)} \end{pmatrix}, \quad M = \begin{pmatrix} 0 \\ -\int_0^{l_1} \left( \frac{1}{\gamma_1(x)} \right) I_{x,l_1}^\beta \omega_1(x) dx \\ \vdots \\ -\int_0^{l_n} \left( \frac{1}{\gamma_n(x)} \right) I_{x,l_n}^\beta \omega_n(x) dx \end{pmatrix}.$$

Indeed, the determinant is

$$\det P = -\sum_{i=1}^{n} \frac{1}{\int_0^{l_i} \frac{dx}{\gamma_i(x)}} \cdot \prod_{j=1}^{n} \int_0^{l_j} \frac{dx}{\gamma_j(x)} \neq 0.$$

So, $\mu = P^{-1}M$. We find unknown coefficients in the expression of the solution. This way, the inverse operator is found.

Further, from (3.7) we obtain

$$0 = \int_{t_1}^{t} d\tau \int_{\Gamma} \left( \partial_{0,t}^{\alpha} v + \partial_{-}^{\beta} \left( \gamma D_{+}^{\beta} v \right) \right) \cdot \partial_{0,t}^{\alpha} \left( \partial_{-}^{\beta} \left( \gamma D_{+}^{\beta} v \right) \right) d\Gamma =$$

$$= \int_{t_1}^{t} d\tau \int_{\Gamma} \partial_{0,t}^{\alpha} \left( \partial_{-}^{\beta} \left( \gamma D_{+}^{\beta} v \right) \right) \cdot \partial_{0,t}^{\alpha} v \, d\Gamma + \int_{t_1}^{t} d\tau \int_{\Gamma} \partial_{0,t}^{\alpha} \left( \partial_{-}^{\beta} \left( \gamma D_{+}^{\beta} v \right) \right) \cdot \partial_{-}^{\beta} \left( \gamma D_{+}^{\beta} v \right) d\Gamma \geq$$

$$\geq p_1 \cdot \int_{t_1}^{t} d\tau \int_{\Gamma} \left[ D_{+}^{\beta} \left( \partial_{0,t}^{\alpha} v \right) \right]^2 d\Gamma + \frac{1}{2} \partial_{t_1,t}^{\alpha-1} \int_{\Gamma} \left[ \partial_{-}^{\beta} \left( \gamma D_{+}^{\beta} v \right) \right]^2 d\Gamma.$$

From the last inequality, we get $Lv := \partial_{-}^{\beta} \left( \gamma D_{+}^{\beta} v \right) = 0$ for a. e $t_1 < t < T$, $x \in \Gamma$. Consequently, $\omega = \partial_{t_1,t}^{\alpha} Lv = 0$ for a. e. $t_1 < t < T, x \in \Gamma$. Taking to account, that $t_1$ is arbitrary number in $(0, T)$, we get $\omega = 0$ in $L_2(G_T)$. The lemma is proven.

We can conclude from this that *Problem D* (the direct problem) has a unique generalized solution in $V(G_T)$. Thus, the Theorem 3.1 is proven.

## 4. SOLVABILITY OF THE INVERSE PROBLEM

Let the following conditions be satisfied

(K1) $g(x,t) \in L_{\infty}(0,T; L_2(\Gamma))$, $\operatorname*{ess\,sup}_{0 \leq t \leq T} \|g(\cdot,t)\|_{L_2(\Gamma)} \leq c$,

$\eta(x) \in H_+^{\beta}(\Gamma)$, $I_{0,x}^{1-\beta} \eta_k(0) = I_{0,x}^{1-\beta} \eta_j(0)$, $k \neq j$, $k, j = \overline{1,n}$, $I_{0,x}^{1-\beta} \eta_k(l_k) = 0$, $k = \overline{1,n}$,

$\|D_{0,x}^{\beta} \eta(x)\|_{L_2(\Gamma)} = m > 0$, $\psi(t) \in H^{\alpha}(0,T)$, $|g^*(t)| \geq q > 0$,

where

$$g^*(t) = \sum_{k=1}^{n} \int_0^{l_k} \eta_k(x) g_k(x,t) dx, \quad t \in (0,T].$$

Multiplying both sides of equation (1.7) by the function $\eta(x)$ and integrating over $\Gamma$, we obtain

$$\sum_{k=1}^{n} \int_0^{l_k} \eta_k(x) \partial_{0,t}^{\alpha} z_k(x,t) dx + \sum_{k=1}^{n} \int_0^{l_k} \gamma_k(x) D_{0,x}^{\beta} z_k(x,t) D_{0,x}^{\beta} \eta_k(x) dx = f(t) g^*(t).$$

Hence, we can find $f$ from this

$$f = \frac{1}{g^*(t)} \sum_{k=1}^{n} \int_0^{l_k} \eta_k(x) \partial_{0,t}^{\alpha} z_k(x,t) dx + \frac{1}{g^*(t)} \sum_{k=1}^{n} \int_0^{l_k} \gamma_k(x) D_{0,x}^{\beta} z_k(x,t) D_{0,x}^{\beta} \eta_k(x) dx.$$

We rewrite the last equation in the following form

$$f = Bf + \frac{\partial_{0,t}^{\alpha} E(t)}{g^*(t)}, \qquad (4.1)$$

where $E(t)$ is defined in (1.8),

$$(Bf)(t) = \frac{1}{g^*}\left\{\sum_{k=1}^{n}\int_0^{l_k}\gamma_k(x)D_{0,x}^{\beta}z_k(x,t)D_{0,x}^{\beta}\eta_k(x)dx\right\}, \qquad (4.2)$$

$$B: L_2(0,T) \to L_2(0,T).$$

It is easy to see that if conditions (K1) are satisfied then the overdetermination condition (1.8) can be equivalently replaced to the equation (4.1).

**Theorem 4.1.** Let condition (K1) hold. If $h(x,t) \in L_2(G_T)$, then the inverse problem (1.1)–(1.6) has a unique generalized solution $\{u,f\} \in V \times L_2(0,T)$.

**Proof.** As we proved the unique solvability of the corresponding direct problem, without loss of generality, we put $h=0$, i.e., we consider *Problem I*. First, we show that the resolvent operator $(I-B)^{-1}$ is bounded and continuous as a map from $L_2(0,T)$ to $L_2(0,T)$.

We multiply each of equations (1.7) by corresponding $\partial_{0,t}^{\alpha}z_k$ and integrate over $G_t$ to get

$$\int_0^t d\tau \int_\Gamma \left(\partial_{0,t}^{\alpha}z(x,t)\right)^2 d\Gamma + \int_0^t d\tau \int_\Gamma \partial_{-}^{\beta}\left(\gamma(x)D_{+}^{\beta}z(x,t)\right)\partial_{0,t}^{\alpha}z(x,t)d\Gamma =$$

$$= \int_0^t f(\tau)d\tau \int_\Gamma g(x,t)\partial_{0,t}^{\alpha}z(x,t)d\Gamma.$$

Taking into account the inequality (3.5), we can write as

$$\left\|\partial_{0,t}^{\alpha}z\right\|_{L_2(G_t)}^2 + \frac{p_1}{2}\partial_{0,t}^{\alpha-1}\left\|D_{+}^{\beta}z\right\|_{L_2(\Gamma)}^2 \le \frac{1}{2}\int_0^t f^2(\tau)d\tau \int_\Gamma g^2(x,\tau)d\Gamma + \frac{1}{2}\left\|\partial_{0,t}^{\alpha}z\right\|_{L_2(G_t)}^2.$$

Then we get

$$\left\|\partial_{0,t}^{\alpha}z\right\|_{L_2(G_t)}^2 + p_1 \cdot \partial_{0,t}^{\alpha}\left\|D_{+}^{\beta}z\right\|_{L_2(G_t)}^2 \le c^2\left\|f(t)\right\|_{L_2(0,t)}^2.$$

From the last inequality, we have

$$\left\|D_{+}^{\beta}z\right\|_{L_2(G_t)}^2 \le \frac{c^2}{p_1}I_{0,t}^{\alpha}\left\|f(t)\right\|_{L_2(0,t)}^2. \qquad (4.3)$$

From (4.3) and taking into account (K1) we find

$$\left\|(Bf)(t)\right\|_{L_2(0,t)}^2 \le \frac{p_2^2 m^2}{q^2}\left\|D_{+}^{\beta}z\right\|_{L_2(0,t)}^2 \le CI_{0,t}^{\alpha}\left\|f\right\|_{L_2(0,t)}^2, \qquad (4.4)$$

where $C = \dfrac{p_2^2 m^2 c^2}{p_1 q^2}$.

Now, iterating the inequality (4.4) $j$ times we get

$$\left\|(B^j f)(t)\right\|_{L_2(0,t)}^2 \le C^j I_{0,t}^{\alpha}\left\|B^{j-1}f\right\|_{L_2(0,t)}^2 = C^j I_{0,t}^{j\alpha}\left\|f\right\|_{L_2(0,t)}^2, \quad j=1,2,3,\ldots$$

From the last inequality, taking into account that the function $\tilde{f}(t) = \left\|f\right\|_{L_2(0,t)}^2$ is a nonnegative and non-decreasing function on $t \in [0,T]$, we have

$$\left\|(B^j f)(t)\right\|_{L_2(0,t)}^2 \le C^j \left\|f\right\|_{L_2(0,t)}^2 \cdot I_{0,t}^{j\alpha}1 = \frac{C^j t^{j\alpha}}{\Gamma(j\alpha+1)}\left\|f\right\|_{L_2(0,t)}^2.$$

Accordingly, we have

$$\left\|(I-B)^{-1}f\right\|_{L_2(0,T)} \leq \sum_{j=0}^{+\infty} \frac{\left(\sqrt{CT^\alpha}\right)^j}{\sqrt{\Gamma(j\alpha+1)}}\|f\|_{L_2(0,t)}.$$

So, it follows that the resolvent operator $(I-B)^{-1}: L_2(0,T) \to L_2(0,T)$ is bounded and continuous mapping and

$$f = (I-B)^{-1}\left[\frac{\partial_{0,t}^\alpha E(t)}{g^*(t)}\right].$$

Now, as $f(t)g(x,t) \in L_2(G_T)$, according to Theorem 1 we have $u(x,t) \in V(G_t)$. This completes the proof of the Theorem 2.

## 5. CONCLUSION

In this work, we considered the direct problem and inverse source problem with integral overdetermination condition for the space-time fractional parabolic equation given on a metric star graph. The existence and uniqueness of a generalized solution to the direct problem is proven by the operator method and using a priori estimates. The inverse problem also reduced another operator equation. It is shown that the resolvent of this operator is well-defined in the inverse source problem given by the overdetermination condition.

## CONFLICT OF INTEREST

The authors declare that they have no conflicts of interest.

FIGURE CAPTIONS

**Fig. 1.** Metric star graph.

FIGURES

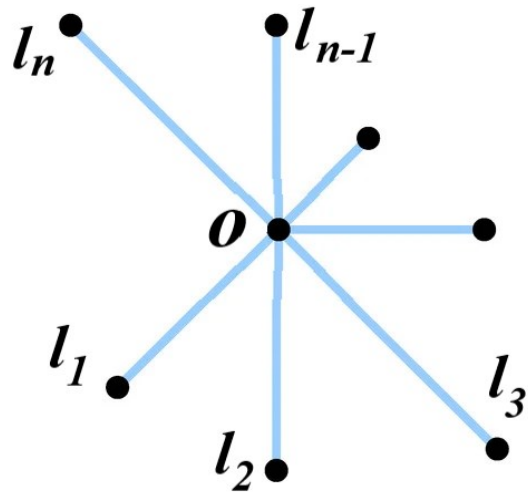

Fig. 1.